\documentclass[10pt,a4paper]{article}
\usepackage{setspace}
\usepackage{amsfonts,amsmath,amssymb}
\usepackage{graphicx}
\usepackage{sectsty}
\usepackage{fullpage}
\newtheorem{theorem}{Theorem}[section]
\newtheorem{lemma}[theorem]{Lemma}

\newtheorem{proposition}[theorem]{Proposition}
\newtheorem{corollary}[theorem]{Corollary}
\newenvironment{proof}[1][Proof]{\begin{trivlist}
\item[\hskip \labelsep {\bfseries #1}]}{\end{trivlist}}
\numberwithin{equation}{section}

\sectionfont{\large}

\begin{document}
\title{Uniqueness on the Class of Odd-Dimensional Starlike Obstacles with Cross Section Data }
\author{Lung-Hui Chen$^1$}\maketitle\footnotetext[1]{Department of
Mathematics, National Chung Cheng University, 168 University Rd.,
Min-Hsiung, Chia-Yi County 621, Taiwan. Email:
mr.lunghuichen@gmail.com.  Fax: 886-5-2720497. The author is
supported by NSC Grant 97-2115-M194-010-MY2.}
\begin{abstract}
We determine the uniqueness on starlike obstacles by using the
cross section data. We see cross section data as spectral measure
in polar coordinate at far field. Cross section scattering data
suffice to give the local behavior of the wave trace. These local
trace formulas contain the geometric information on the obstacle.
Local wave trace behavior is connected to the cross section
scattering data by Lax-Phillips' formula. Once the scattering data
are identical from two different obstacles, the short time
behavior of the localized wave trace is expected to give identical
heat/wave invariants.
\end{abstract}

\section{Introduction and the Statement of Main Result}
Let $H$ be an embedded hypersurface in $\mathbb{R}^n$ such that
\begin{equation}
\mathbb{R}^n\setminus H=\Omega\cup\mathcal{O},\mbox{ with
}\overline{\mathcal{O}}\mbox{ compact and }\overline{\Omega}\mbox{
connected },
\end{equation}where both $\mathcal{O}$ and $\Omega$ are open. We
call $\mathcal{O}$ an obstacle and $\Omega$ as its exterior.
Without loss of generality, we assume $\mathcal{O}$ contains the
origin.

\par
Mathematically, exterior scattering problem is formulated as
follows. Let $u\in\mathcal{C}^\infty(\Omega)$ be the solution to
the following exterior problem
\begin{eqnarray}\label{dirichlet}
\left\{
\begin{array}{cc}
    \Delta u+\lambda^2u=0 & \hbox{ in }\Omega,\\
    u=0& \hbox{ on }H,
  \end{array}\right.
\end{eqnarray}
for Dirichlet condition;
\begin{eqnarray}\label{neumann} \left\{
\begin{array}{cc}
    \Delta u+\lambda^2u=0 & \hbox{ in }\Omega,\\
    \frac{\partial}{\partial \nu}u=0& \hbox{ on }H,
  \end{array}\right.
\end{eqnarray}
for Neumann condition. Let us call the Laplacian defined above as
$\Delta_\mathcal{O}$.
\par
Let $u=u(x,\omega,\lambda)$ be the corresponding incoming
solution. Let $x:=|x|\frac{x}{|x|}=r\theta$. We have the following
asymptotic behavior when $r:=|x|\rightarrow\infty$, for $\lambda$
near $\mathbb{R}$,
\begin{equation}\label{1.4}
u(x,\omega,\lambda)=e^{-i\lambda\omega\cdot x}+v_\omega(x),
\end{equation}
where
\begin{equation}\label{1.5}
v_\omega(x):=\frac{e^{i\lambda
r}}{r^{(n-1)/2}}(A(\lambda,\theta,\omega)+O(\frac{1}{r})), \mbox{
as } r=|x|\rightarrow\infty.
\end{equation}
The function
$A(\lambda,\theta,\omega)\in\mathcal{C}^\infty(\mathbb{R}\setminus\{0\}\times
\mathbb{S}^{n-1}\times \mathbb{S}^{n-1})$ is the scattering
amplitude related to obstacle $\mathcal{O}$. Note that, in the
sense of distribution on $\mathbb{S}^{n-1}_\theta$,
\begin{equation}\label{1.6}
e^{-i\lambda r\theta\cdot\omega}=(\frac{2\pi
i}{\lambda})^{\frac{n-1}{2}}r^{-\frac{n-1}{2}}\{e^{-i\lambda r}
(\delta_\omega(\theta)+O(\frac{1}{r}))+i^{n-1}e^{i\lambda r}
(\delta_{-\omega}(\theta)+O(\frac{1}{r}))\}.
\end{equation}
We define scattering matrix $S(\lambda)$ as the operator with
$\mathcal{C}^\infty$-Schwartz kernel
\begin{equation}\label{1.7}
S(\lambda,\omega,\theta)=\delta_\omega(\theta)+c_n
\lambda^{(n-1)/2}\overline{A(\lambda,-\theta,\omega)},\hspace{2pt}
c_n=(2\pi)^{-\frac{n-1}{2}}e^{-\frac{i\pi}{4}(n-1)}.
\end{equation}
A scattering matrix in this form is close to the one in Lax and
Phillips \cite{Lax and Phillips}. In Melrose \cite[p.23]{Melrose},
we have the "absolute scattering matrix" defined as
\begin{equation}\label{18}
\widehat{S}(\lambda,\omega,\theta)=i^{n-1}\delta_{-\omega}(\theta)+\overline{c_n}
\lambda^{(n-1)/2}A(-\lambda,\theta,\omega),\mbox{ where
}A(\lambda,\theta,\omega)\in\mathcal{C}^\infty(\mathbb{R}\setminus\{0\}\times
\mathbb{S}^{n-1}\times \mathbb{S}^{n-1}),
\end{equation}
which is obtained by comparing the coefficient of the
$e^{-i\lambda r}$ and the one of $e^{i\lambda r}$ as an operator
image in~(\ref{1.4}).
\par
Alternatively, scattering matrix can be derived from Poisson
operator $P(\lambda):\mathcal{L}^2(\mathbb{S}^{n-1})\rightarrow
\mathcal{L}^2(\Omega)$, which has $\mathcal{L}^2$-kernel defined
as
\begin{equation}\label{1.8}
P(\lambda,x,\omega):=
\lambda^{\frac{n-1}{2}}c_nu(x,\omega,\lambda).
\end{equation}

\par
To understand $S(\lambda)$, we begin with the spectral theory of
its resolvent. We define
$$\mathcal{P}:=\{\lambda\in \mathbb{C} |\Im\lambda>0\}$$ as the
physical plane in this paper. According to Sj\"{o}strand and
Zworski \cite{Sjostand and Zworski}, the scattering matrices
$S(\lambda)$ has the meromorphic extension to $\mathbb{C}$ when
$n$ is odd; $\Lambda$, logarithmic plane, when $n$ is even. We use
\begin{equation}
(\Delta_\mathcal{O}-\lambda^2)^{-1}:\mathcal{L}^2(\Omega)\rightarrow
\mathcal{H}^2(\Omega)
\end{equation}
as the scattered resolvent, imposed with scatterers described
above, which is defined over $\mathcal{P}$ by spectral analysis.
As a special case of black-box formalism of Sj\"{o}strand and
Zworski \cite{Sjostand and Zworski},
$(\Delta_\mathcal{O}-\lambda^2)^{-1}:\mathcal{L}^2(\Omega)\rightarrow
\mathcal{H}^2(\Omega)$ meromorphically extends from $\mathcal{P}$,
$\lambda^2\not\in Spec_{pp}(\Delta_\mathcal{O})$, to $\mathbb{C}$
if n is odd; to $\Lambda$, the logarithmic plane, if n is even, as
an operator
\begin{equation}
R(\lambda):\mathcal{L}^2_{\rm comp}(\Omega)
\rightarrow\mathcal{H}^2_{\rm loc}(\Omega)
\end{equation}
In this paper, $n$ is odd. $R(\lambda)$ shares the same spectral
structure as the corresponding scattering matrix $S(\lambda)$. The
resolvent operator $R(\lambda)$ that we will use in this paper are
meromorphically extended. That means they are spatially cut offs.
So do the wave groups.
\par
Inverse scattering theory asks for the information on the
scatterer $\mathcal{O}$ given the knowledge provided by
$S(\lambda)$. In particular, let $\mathcal{O}^1$ and
$\mathcal{O}^2$ be two obstacles, uniqueness problem asks that if
$\mathcal{O}^1=\mathcal{O}^2$ given the information of
$S^1(\lambda,\omega,\theta)=S^2(\lambda,\omega,\theta)$ on partial
or complete set of $(\lambda,\omega,\theta)\in
\mathbb{F}\times\mathbb{S}^{n-1}\times\mathbb{S}^{n-1}$, where
$\mathbb{F}=\mathbb{C}$ or $\Lambda$. Theoretically, the
singularity structure of the scattering matrix may determine the
obstacle. We refer to Isakov's papers \cite{Isakov,Isakov2} for an
earlier review on the uniqueness and the stability results for
obstacle scattering. We refer the inverse scattering problem for
convex bodies to \cite[Theorem 3.2]{Colton and Sleeman} in which
the case for sound-hard and convex obstacle are discussed.
However, there are not too many results on the inverse scattering
problem by cross section data. There are some numerical results
involved with the determination of the obstacle $\mathcal{O}$ by
the corresponding scattering cross section which is defined in
this paper as
\begin{equation}
C(\lambda,\theta):=\int_{\mathbb{S}^{n-1}}|A(\lambda,\theta,\theta')|^2d\theta'.\label{cs}
\end{equation}
As asked by Colton and Sleeman \cite{Colton and Sleeman}, how far
can we determine the obstacle $\mathcal{O}$ from the cross section
$C(\lambda,\theta)$ provided the obstacle is convex and
sound-soft? In \cite{Colton and Sleeman}, the capacity of the
obstacle $\mathcal{O}$ and the areas of the shadow projections of
$\mathcal{O}$ of all directions can be uniquely determined. In
this paper, we will connect the cross section $C(\lambda,\theta)$
to spectral measure and Birman-Krein formula or Hille-Yoshida
formula. Therefore, some geometric invariants can be obtained via
the asymptotic spectral expansion of heat/wave propagator in short
time. Cross section $C(\lambda,\theta)$ can be interpreted as a
directional spectral measure propagating along direction $\theta$.
 How far can we go to tell the geometric difference of these two
obstacles by comparing their heat/wave invariants?
\par
Assuming the boundary defining function of obstacle
$\mathcal{O}^k$, denoted as $x^k$, $k=1,2$, is of the form
\begin{equation}
x^k:=r^k(\theta),\mbox{ where }\theta\in \mathbb{S}^{n-1}\mbox{
and }r^k\in \mathbb{R}^+,
\end{equation}
we state the main result in this paper as
\begin{theorem}
\label{1.1} Let $\mathcal{O}^k$, $k=1,2$, be two starlike
obstacles containing the origin in $\mathbb{R}^n$, $n\geq3$, odd,
with smooth boundary imposed either with~(\ref{dirichlet})
or~(\ref{neumann}). Let $C^k(\lambda,\theta)$ be the cross section
data corresponding to obstacle $\mathcal{O}^k$. If, in the
neighborhood of one fixed $\lambda_0\in\mathbb{R}\setminus\{0\}$,
\begin{equation}
C^1(\lambda,\theta)=C^2(\lambda,\theta),\hspace{2pt}\mbox{for all
}\theta\in\mathbb{S}^{n-1},\label{assumption}
\end{equation}
then we have
\begin{equation}
r^1(\theta)=r^2(\theta), \forall\theta\in\mathbb{S}^{n-1}.
\end{equation}
In particular,
\begin{equation}
\mathcal{O}^1=\mathcal{O}^2.
\end{equation}
\end{theorem}

\section{On the Both Sides of Trace Formulas}

We recall the following theorem from Lax-Phillips' scattering
theory \cite{Lax and Phillips}.
\begin{theorem}[Lax and Phillips]
The scattering amplitude satisfies the following relations:
\begin{eqnarray}\nonumber
&&A(\lambda,-\omega,\theta)+(-1)^{\frac{n-1}{2}}\overline{A(\lambda,-\theta,\omega)}\\\nonumber
&=&-(\frac{\lambda}{2\pi
i})^{\frac{n-1}{2}}\int_{\mathbb{S}^{n-1}}
A(\lambda,-\omega,\theta')\overline{A(\lambda,-\theta,\theta')}d\theta'\\
&=&-(\frac{\lambda}{2\pi
i})^{\frac{n-1}{2}}\int_{\mathbb{S}^{n-1}}
\overline{A(\lambda,-\theta',\omega)}A(\lambda,-\theta',\theta)d\theta'.\label{1.16}
\end{eqnarray}
Moreover,
\begin{equation}
A(\lambda,\omega,\theta)=A(\lambda,\theta,\omega)\mbox{  and
}A(-\lambda,\omega,\theta)=\overline{A(\lambda,\omega,\theta)}.\label{L2}
\end{equation}
\end{theorem}
Furthermore,
\begin{lemma}\label{22}
In the setting from~(\ref{1.4}) to~(\ref{1.7}), we have
\begin{equation}
S^{-1}(\lambda,\omega,\theta)=S(-\lambda,\omega,\theta)
\end{equation}
and, formally as coefficient in the spectral
expansion~(\ref{1.4}),
\begin{equation}
\widehat{S}^{-1}(\lambda,\omega,\theta)=
i^{n-1}\delta_{-\omega}(\theta) +(\frac{\lambda}{2\pi
i})^{\frac{n-1}{2}}e^{-2i\lambda r}A(-\lambda,\theta,-\omega).
\end{equation}
\end{lemma}
\begin{proof}
The first equality comes from Shenk and Thoe \cite{Shenk}. The
proof on the second equality is a straightforward inverse
correspondence in~(\ref{1.4}): let $d_n:=c_n^{-1}$. We
alternatively rewrite~(\ref{1.4}) and~(\ref{1.5}),
\begin{eqnarray}\nonumber
u(x,\omega,\lambda)&=&
\lambda^{-\frac{n-1}{2}}d_nr^{-\frac{n-1}{2}}e^{-i\lambda
r}\delta_\omega(\theta)+\lambda^{-\frac{n-1}{2}}d_ni^{n-1}r^{-\frac{n-1}{2}}e^{i\lambda
r}\delta_{-\omega}(\theta)\\&&+r^{-\frac{n-1}{2}}e^{-i\lambda
r}e^{2i\lambda r}A(\lambda,\theta,\omega)+\cdots.
\end{eqnarray}
Hence,
\begin{eqnarray}\nonumber
u(x,\omega,\lambda)&=&
i^{n-1}r^{-\frac{n-1}{2}}\lambda^{-\frac{n-1}{2}}d_n\{e^{i\lambda
r}\delta_{-\omega}(\theta)\\&&+e^{-i\lambda
r}[i^{n-1}\delta_\omega(\theta)+i^{n-1}\lambda^{\frac{n-1}{2}}d_n^{-1}e^{2i\lambda
r}A(\lambda,\theta,\omega)]\}+\cdots.
\end{eqnarray}
Observing the correspondence from the coefficient of $e^{i\lambda
r}$ to $e^{-i\lambda r}$, we obtain the kernel of
$\widehat{S}^{-1}(\lambda)$. $\Box$
\end{proof}

To connect $C(\lambda,\theta)$ to spectral analysis, we look at
the following lemma.
\begin{lemma}\label{11}
Under the theorem assumption, we let $R^k(\lambda,x,y)$ be the
resolvent kernel corresponding to $\mathcal{O}^k$ with exterior
$\Omega^k$, $k=1,2$. Then, in a neighborhood of $\lambda_0$ in
$0i+\mathbb{R}$,
\begin{eqnarray}2\lambda
\{R^1(\lambda,x,y)-R^1(-\lambda,x,y)\}= 2\lambda
\{R^2(\lambda,x,y)-R^2(-\lambda,x,y)\},\forall
x,y\in\Omega^1\cap\Omega^2.
\end{eqnarray}
\end{lemma}
\begin{proof}
Starting with
\begin{equation}
[R(\lambda)-R(-\lambda)]d\lambda^2=(2\pi)^{-1}P(\lambda)P^\ast(\lambda)d\lambda,
\end{equation}in which either quantities can serve as the
definition of spectral measure. See Reed and Simon \cite{Reed and
Simon}. Letting $R^k$, $P^k$, $u^k$, $A^k$ and $C^k$ be the
corresponding quantities related to obstacle $\mathcal{O}^k$, we
have
\begin{eqnarray}\nonumber
&&2\lambda\{R^1(\lambda,x,x)-R^1(-\lambda,x,x)\}-2\lambda\{R^2(\lambda,x,x)-R^2(-\lambda,x,x)\}\\\nonumber
&=&\frac{1}{2\pi}\int_{\mathbb{S}^{n-1}}P^1(\lambda,x,\omega)\overline{P^1(\lambda,\omega,x)}d\omega\nonumber
-\{\mbox{similar terms from } R^2\}\\\nonumber
&=&\frac{\lambda^{n-1}|c_n|^{2}}{2\pi}\int_{\mathbb{S}^{n-1}}
1+\frac{e^{-i\lambda\omega\cdot x }e^{-i\lambda
r}}{r^{\frac{n-1}{2}}}\overline{A^1(\lambda,\omega,\theta)}+\frac{e^{i\lambda\omega\cdot
x }e^{i\lambda r}}{r^{\frac{n-1}{2}}}A^1(\lambda,\theta,\omega)
+\frac{\overline{A^1(\lambda,\omega,\theta)}A^1(\lambda,\omega,\theta)}{r^{n-1}}d\omega\nonumber\\\nonumber
&&-\{\mbox{similar terms from } R^2\}+O(\frac{1}{r^{n}}),\mbox{ by
}~(\ref{1.4}),~(\ref{1.5}),~(\ref{1.6}),~(\ref{1.8}),\\\nonumber
&=&\frac{\lambda^{n-1}|c_n|^{2}}{2\pi}\int_{\mathbb{S}^{n-1}}\frac{-e^{i\lambda\omega\cdot
x}e^{-i\lambda
r}}{r^{\frac{n-1}{2}}}\overline{A^1(\lambda,\theta,\omega)}+\frac{e^{i\lambda\omega\cdot
x}e^{i\lambda r}}{r^{\frac{n-1}{2}}}A^1(\lambda,\theta,\omega)
+\frac{\overline{A^1(\lambda,\theta,\omega)}A^1(\lambda,\theta,\omega)}{r^{n-1}}d\omega\nonumber\\
&&-\{\mbox{similar terms from }
R^2\}+O(\frac{1}{r^{n}})\nonumber\\
&=&\frac{\lambda^{n-1}|c_n|^{2}}{2\pi}\frac{e^{-i\lambda
r}}{r^{\frac{n-1}{2}}}\int_{\mathbb{S}^{n-1}}e^{-i\lambda\omega\cdot
x}\overline{A^1(\lambda,\theta,\omega)}d\omega+\frac{\lambda^{n-1}|c_n|^{2}}{2\pi}\frac{e^{i\lambda
r}}{r^{\frac{n-1}{2}}}\int_{\mathbb{S}^{n-1}}e^{i\lambda\omega\cdot
x}A^1(\lambda,\theta,\omega)d\omega\nonumber\\&&
+\frac{(2\pi)^{-1}\lambda^{n-1}|c_n|^{2}}{r^{n-1}}\int_{\mathbb{S}^{n-1}}
\overline{A^1(\lambda,\theta,\omega)}A^1(\lambda,\theta,\omega)d\omega\nonumber\\
&&-\{\mbox{similar terms from } R^2\}+O(\frac{1}{r^{n}})
.\label{124}
\end{eqnarray}
We compute this term by term. Using~(\ref{1.6}),
\begin{eqnarray}\nonumber
&&\int_{\mathbb{S}^{n-1}} e^{i\lambda\omega\cdot x}A^k(\lambda,\theta,\omega)d\omega\\
&\underset{x\rightarrow\infty}{\rightarrow} &
\int_{\mathbb{S}^{n-1}}
(\frac{2\pi}{i\lambda})^{\frac{n-1}{2}}r^{-\frac{n-1}{2}}
[e^{i\lambda r}\delta_\omega(\theta)+i^{n-1}e^{-i\lambda
r}\delta_{-\omega}(\theta)]
A^k(\lambda,\theta,\omega)d\omega\nonumber\\
&=&(\frac{2\pi }{i\lambda
})^{\frac{n-1}{2}}r^{-\frac{n-1}{2}}e^{i\lambda
r}A^k(\lambda,\theta,\theta)+(\frac{2\pi }{i\lambda
})^{\frac{n-1}{2}}r^{-\frac{n-1}{2}}i^{n-1}e^{-i\lambda
r}A^k(\lambda,\theta,-\theta)+\cdots.
\end{eqnarray}
Taking conjugate,
\begin{eqnarray}\nonumber
&&\int_{\mathbb{S}^{n-1}} e^{-i\lambda\omega\cdot x}\overline{A^k(\lambda,\theta,\omega)}d\omega\\
&\underset{x\rightarrow\infty}{\rightarrow}& (\frac{2\pi
i}{\lambda })^{\frac{n-1}{2}}r^{-\frac{n-1}{2}}e^{-i\lambda
r}\overline{A^k(\lambda,\theta,\theta)}+(\frac{2\pi
i}{\lambda})^{\frac{n-1}{2}}r^{-\frac{n-1}{2}}i^{n-1}e^{i\lambda
r}\overline{A^k(\lambda,\theta,-\theta)}+\cdots.
\end{eqnarray}
However, from the identities in Lemma~\ref{22},
\begin{eqnarray}
&&S^k(-\lambda,\omega,\theta)= \delta_{\omega}(\theta)
+(\frac{-\lambda}{2\pi i})^{\frac{n-1}{2}}A^k(\lambda,-\theta,\omega);\\
&&\{\widehat{S}^k\}^{-1}(\lambda,\omega,\theta)=
i^{n-1}\delta_{-\omega}(\theta) +(\frac{\lambda}{2\pi
i})^{\frac{n-1}{2}}e^{-2i\lambda
r}A^k(-\lambda,\theta,-\omega).\label{2.15}
\end{eqnarray}
Using~(\ref{18}) we have
\begin{equation}
\overline{A^k(\lambda,\theta,\omega)}= e^{2i\lambda
r}A^k(\lambda,\theta,-\omega).
\end{equation}
Let $\omega=-\theta$. We obtain
\begin{equation}\label{2222}
\overline{A^k(\lambda,\theta,-\theta)}= e^{2i\lambda
r}A^k(\lambda,\theta,\theta).
\end{equation}

Therefore, as $|x|\rightarrow\infty$,
\begin{eqnarray}\nonumber
&&\lambda^{n-1}|c_n|^{2}\int_{\mathbb{S}^{n-1}}\frac{e^{i\lambda\omega\cdot
x}e^{i\lambda
r}}{r^{\frac{n-1}{2}}}A^k(\lambda,\theta,\omega)+\frac{e^{-i\lambda\omega\cdot
x}e^{-i\lambda
r}}{r^{\frac{n-1}{2}}}\overline{A^k(\lambda,\theta,\omega)}d\omega\\\nonumber
&=&\lambda^{\frac{n-1}{2}}\overline{c_n}\frac{1}{r^{n-1}}(A^k(\lambda,\theta,-\theta)
+(-1)^{\frac{n-1}{2}}\overline{A^k(\lambda,\theta,-\theta)})\\&&+\lambda^{\frac{n-1}{2}}
c_n\frac{1}{r^{n-1}}(e^{2i\lambda
r}A^k(\lambda,\theta,\theta)+(-1)^{\frac{n-1}{2}}e^{-2i\lambda
r}\overline{A^k(\lambda,\theta,\theta)})+O(\frac{1}{r^{n}})\nonumber\mbox{,
using~(\ref{1.16}),~(\ref{2222}),}\\
&=&-2(\frac{\lambda}{2\pi})^{n-1}r^{-(n-1)}C^k(\lambda,\theta)+O(\frac{1}{r^{n}}).\label{1.28}
\end{eqnarray}
Hence,~(\ref{124}) and~(\ref{1.28}) sum up to give
\begin{eqnarray}\nonumber
&&2\lambda\{R^1(\lambda,x,x)-R^1(-\lambda,x,x)\}-2\lambda\{R^2(\lambda,x,x)-R^2(-\lambda,x,x)\}\\
&\underset{x\rightarrow\infty}{\rightarrow}&(\frac{-
1}{2\pi})r^{-(n-1)}(\frac{\lambda}{2\pi})^{n-1}
\{C^1(\lambda,\theta)-C^2(\lambda,\theta)\}+O(\frac{1}{r^{n}}).\label{1.23}
\end{eqnarray}
\par
Furthermore, we see that
$\int_{|x|=s}2\lambda\{R^1(\lambda,x,y)-R^1(-\lambda,x,y)-R^2(\lambda,x,y)+R^2(-\lambda,x,y)\}dS_x$
is a solution of the exterior problem~(\ref{dirichlet})
or~(\ref{neumann}) for $|y|\gg d$,  $\forall s\in[c,d]\subset
\mathbb{R},c\gg1$. Therefore, using Jensen's inequality, for some
constant $C$ depending only on $n$ and $s$,
\begin{eqnarray}\nonumber
&&\int_{|y|=r}\{\int_{|x|=s}2\lambda\{R^1(\lambda,x,y)-R^1(-\lambda,x,y)-R^2(\lambda,x,y)+R^2(-\lambda,x,y)\}dS_x\}^2dS_y\\\nonumber
&\leq&
C\int_{|y|=r}\int_{|x|=s}\{2\lambda\{R^1(\lambda,x,y)-R^1(-\lambda,x,y)-R^2(\lambda,x,y)+R^2(-\lambda,x,y)\}\}^2dS_x dS_y\\
&\leq&C\int_{|y|=r}\int_{|x|=r}|2\lambda\{R^1(\lambda,x,y)-R^1(-\lambda,x,y)-R^2(\lambda,x,y)+R^2(-\lambda,x,y)\}|^2dS_xdS_y\nonumber\\
&\leq&C\{\int_{|x|=r}|2\lambda\{R^1(\lambda,x,x)-R^1(-\lambda,x,x)-R^2(\lambda,x,x)+R^2(-\lambda,x,x)\}|dS_x\}^2,\label{2.18}
\end{eqnarray}
where the last inequality comes from the fact that Hilbert-Schmidt
norm is controlled by trace norm. The theorem assumption
$C^1(\lambda,\theta)=C^2(\lambda,\theta)$, ~(\ref{1.23})
and~(\ref{2.18}) yield
\begin{eqnarray}\nonumber
&&\int_{|y|=r}\{\int_{|x|=s}2\lambda\{R^1(\lambda,x,y)-R^1(-\lambda,x,y)-R^2(\lambda,x,y)+R^2(-\lambda,x,y)\}dS_x\}^2dS_y\\
&\lesssim& C\{\int_{|x|=r}|x|^{-n}dS_x\}^2\nonumber\\
&\leq& Cr^{-2},\mbox{ where }C's\mbox{ are constants}.
\end{eqnarray}
Using the Kato's uniqueness theorem as in Shenk and Thoe
\cite[Lemma 4.4]{Shenk2}, we have
\begin{equation}
\int_{|x|=s}2\lambda\{R^1(\lambda,x,y)-R^1(-\lambda,x,y)-R^2(\lambda,x,y)+R^2(-\lambda,x,y)\}dS_x\equiv0,\forall
y\in\Omega^1\cap\Omega^2,\forall s\in[c,d].
\end{equation}
By Lebesgue integration theory, $\forall
y\in\Omega^1\cap\Omega^2$,
$\{R^1(\lambda,x,y)-R^1(-\lambda,x,y)-R^2(\lambda,x,y)+R^2(-\lambda,x,y)\}$
is $a.e.$ zero with respect to $|x|\in[c,d]$. Since it is
continuous to $x$, provided by Theorem 5.1 part(3) in
\cite{Shenk2}, it is identically zero with respect to
$|x|\in[c,d]$, $\forall y\in\Omega^1\cap\Omega^2$. Now we apply
the unique continuation property of elliptic differential equation
with analytic coefficients. See, Bers, John and Schechter
\cite{Bers}. Here, we have Helmholtz equation as a special case.

\par
Again, using the unique continuation property of Helmholtz
equation with respect to $x$, we have
\begin{equation}
2\lambda\{R^1(\lambda,x,y)-R^1(-\lambda,x,y)-R^2(\lambda,x,y)+R^2(-\lambda,x,y)\}\equiv0,\forall
x,y\in \Omega^1\cap\Omega^2.
\end{equation}
$\Box$
\end{proof}
As a result of continuation outside all possible poles, we have in
particular that
\begin{corollary}\label{2.4}
In $0i+\mathbb{R}$, there exist a cutoff function
$\chi\in\mathcal{C}^\infty_0(\mathbb{R}^n;[0,1])$ which is $1$
near $\mathcal{O}^1\cup\mathcal{O}^2$ such that
\begin{eqnarray}
(1-\chi)\{R^1(\lambda,\cdot,\cdot)-R^1(-\lambda,\cdot,\cdot)\}\equiv
(1-\chi)\{R^2(\lambda,\cdot,\cdot)-R^2(-\lambda,\cdot,\cdot)\}.
\end{eqnarray}
\end{corollary}

Let the naturally regularized wave trace
\begin{equation}
u(t):=2\{\cos{t\sqrt{\Delta_\mathcal{O}}}-\cos{t}\sqrt{\Delta_0}\},\label{2.2}
\end{equation}
where
$\Delta_0:=-\Delta_{\mathbb{R}^n}=-\frac{\partial^2}{\partial
x_1^2}-\frac{\partial^2}{\partial
x_2^2}-\cdots-\frac{\partial^2}{\partial x_n^2}$.
$\cos{t\sqrt{\Delta_\mathcal{O}}}$ has a kernel satisfying the
following Cauchy problem:
\begin{equation}
\left\{%
\begin{array}{ll}
    (\frac{\partial^2}{\partial t^2}+\Delta_\mathcal{O})\cos{t\sqrt{\Delta_\mathcal{O}}}(x,y)=0; \\
    \cos{t\sqrt{\Delta_\mathcal{O}}}(x,y)|_{t=0}=\delta(x-y); \\
    \frac{\partial \cos{t\sqrt{\Delta_\mathcal{O}}}(x,y)}{\partial t}|_{t=0}=0.  \\
\end{array}%
\right.
\end{equation}
Furthermore, $u(t)$ has a distributional trace. See Zworski
\cite{Zworski3}. We recall from Petkov and Stoyanov \cite{Petkov}
that, for a non-trapping obstacle,
\begin{equation}
\mbox{singular support of }Tr\{u(t)\}=\{0\}.\label{2.25}
\end{equation}
Moreover, $\cos{t\sqrt{\Delta_\mathcal{O}}}(x,y)$, $t\geq0$, is
interpreted as the data given at $(0,y)$ received at $(t,x)$ along
the geodesic. Hence, we see
$\cos{t\sqrt{\Delta_\mathcal{O}}}(x,x)$, $t\geq0$, as the data
given at $(-t/2,x)$ received at $(t/2,x)$ along the geodesic.

Furthermore, $u(t)$ has a spectral representation.
\begin{eqnarray}
u(t)&\hspace{-2pt}=\hspace{-2pt}&\int_\mathbb{R}e^{-i\lambda
t}\{R(\lambda)-R^0(\lambda)-R(-\lambda)+R^0(-\lambda)\}d\lambda^2.\label{2.26}
\end{eqnarray}
$\{R(\lambda)-R^0(\lambda)-R(-\lambda)+R^0(-\lambda)\}d\lambda^2$
is the spectral measure, where
$R^0(\lambda):=(\Delta_0-\lambda^2)^{-1}$. There is no Neumann or
Dirichlet eigenvalue of $\Delta_\mathcal{O}$ in obstacle
scattering problem. Furthermore, when $n\geq3$, $0$ is neither a
resonance nor an eigenvalue of $R(\lambda)$. The continuous
spectrum is actually where the scattering phenomenon happens. We
consider the Fourier inversion formula of~(\ref{2.26}) over
$\overline{\mathcal{P}}$
\begin{eqnarray}\label{2.27}
\int_\mathbb{R}e^{i\lambda
t}u(t)dt&\hspace{-2pt}=\hspace{-2pt}&2\lambda
\{R(\lambda)-R^0(\lambda)-R(-\lambda)+R^0(-\lambda)\}.
\end{eqnarray}
Or, locally,
\begin{eqnarray}\label{2.9}
\int_\mathbb{R}e^{i\lambda t}f
u(t)dt&\hspace{-2pt}=\hspace{-2pt}&2\lambda
f\{R(\lambda)-R^0(\lambda)-R(-\lambda)+R^0(-\lambda)\}, \mbox{
where }f\in\mathcal{C}^\infty_0(\Omega).
\end{eqnarray}
We will focus at the  behavior of the localized cutoffed
resolvents on the boundary $H$.

\par
Using Birman-Krein type of theory, we see that, for $\lambda\in
\mathbb{R}$,
\begin{equation}
2\lambda
Tr\{R(\lambda)-R^0(\lambda)-R(-\lambda)+R^0(-\lambda)\}=\sigma'(\lambda).\label{222}
\end{equation}
A general treatment in proving the Birman-Krein theorem in black
box formalism setting when $n$ is odd can be found  in Zworski
\cite{Zworski3}. Therefore, we can rewrite~(\ref{2.27}) in a
distributional sense as
\begin{equation}
\sigma'(\lambda)=\int_\mathbb{R}e^{i\lambda t}Tr
\{u(t)\}dt,\lambda\in0+i\mathbb{R}.
\end{equation}
Locally, we can define
\begin{equation}
{\sigma_f}'(\lambda):=2\lambda
Tr\{f(R(\lambda)-R^0(\lambda)-R(-\lambda)+R^0(-\lambda))\},\forall
f\in\mathcal{C}^\infty_0(\mathbb{R}^n).\label{2.24}
\end{equation}
In this notation, we can convert~(\ref{2.9}) to a local formula:
\begin{equation}
{\sigma_f}'(\lambda)=\int_\mathbb{R}e^{i\lambda t}Tr \{fu(t)\}dt.
\end{equation}
Let ${\sigma_f^{k'}}(\lambda)$ be the quantity corresponding to
$\mathcal{O}^k$. Timing $f$ on the distributional resolvent kernel
$R^k(\lambda,x,x)$ and carrying out the trace integration,
Lemma~\ref{11} tells us
\begin{equation}
{\sigma_f^1}'(\lambda)={\sigma_f^2}'(\lambda),\mbox{ in a
neighborhood of }\lambda_0, \forall
f\in\mathcal{C}^\infty_0(\Omega^1\cap\Omega^2) ).\label{2.31}
\end{equation}

\par
Now we prove
\begin{proposition}\label{2.3}
Under the same assumption as in the introduction, the inverse
Fourier transform corresponding to $\mathcal{O}^k$
$\int_\mathbb{R}e^{i\lambda t}Tr\{fu^k(t)\}dt$, which is valid for
$\lambda\in 0i+\mathbb{R}$, depends only on its short time
behavior in the following sense:
\begin{equation}
\int_\mathbb{R}e^{i\lambda
t}Tr\{f(u^1(t)-u^2(t))\}dt=\int_\mathbb{R}e^{i\lambda
t}Tr\{f(u^1(t)-u^2(t))\}\rho_1(t)dt+\mbox{ rapidly decreasing
term},
\end{equation}
whenever $\lambda\in0i+\mathbb{R}$ and where
$\rho_1(t)\in\mathcal{C}^\infty_0(\mathbb{R};[0,1])$ is a cutoff
function supported at $t=0$. Moreover,
$f\in\mathcal{C}^\infty_0(\mathbb{R}^n;[0,1])$ is $1$ near
$\mathcal{O}^1\cup\mathcal{O}^2$.
\end{proposition}
\begin{proof}
We divide the inverse Fourier transform into three time intervals:
\begin{eqnarray}\nonumber
&&\int_{-\infty}^\infty e^{i\lambda
t}Tr\{f(u^1(t)-u^2(t))\}dt\\\nonumber &:=&\int_{-\infty}^\infty
e^{i\lambda
t}Tr\{f(u^1(t)-u^2(t))\}\rho_1(t)dt+\int_{-\infty}^\infty
e^{i\lambda
t}Tr\{f(u^1(t)-u^2(t))\}\rho_2(t)dt\\&&\nonumber+\int_{-\infty}^\infty
e^{i\lambda
t}Tr\{f(u^1(t)-u^2(t))\}\rho_3(t)dt\\
&:=& I_1(\lambda)+I_2(\lambda)+I_3(\lambda),\label{2.12}
\end{eqnarray}
where $\rho_i\in\mathcal{C}^\infty(\mathbb{R};[0,1])$, $i=1,2,3$.
Let $\rho_1,\rho_3\in\mathcal{C}^\infty(\mathbb{R};[0,1])$ be two
cutoff functions such that $\rho_1$ is $1$ with small compact
support at $t=0$ and $\rho_3$ is $1$ near $t=\pm\infty$. We take
$\rho_1(t)+\rho_2(t)+\rho_3(t)=1$. We take $\beta$ such that ${\rm
supp}(\rho_3(t))\subset(-\infty,-\beta)\cup(\beta,\infty)$.
$\beta$ is to be chosen. This is a partition of unity.
\par
Using Paley-Wiener's theorem for $I_1(\lambda)$,
\begin{equation}
|\int_{-\infty}^\infty e^{i\lambda
t}Tr\{f(u^1(t)-u^2(t))\}\rho_1(t)dt|\leq C(1+|\lambda|)^Ne^{h(-\Im
\lambda)},\label{2.13}
\end{equation}for some $N\in\mathbb{N}$ and for some constant $C$. $h$ is the support
function of $Tr\{f(u^1(t)-u^2(t))\}\rho_1(t)$. We just keep
$I_1(\lambda)$. $N$ will be specified by Ivrii's result
\cite{Ivrii}. $I_1(\lambda)$ is holomorphic and well-defined as a
Fourier-Laplace transform.
\par
We also apply Paley-Wiener's theorem to $I_2(\lambda)$.
By~(\ref{2.25}), for each $\beta$,
$Tr\{f(u^1(t)-u^2(t))\}\rho_2(t)$ is a smooth function with
compact support. By construction $\rho_2(t)$ is the union of two
cutoff functions. One, denoted as $\rho_2^+(t)$, is supported on
$\mathbb{R}^+$ while the other one, denoted as $\rho_2^-(t)$,
supported on $\mathbb{R}^-$. For the first one, we choose
$\Im\lambda>0$, the upper half complex plane, for
\begin{equation}\label{2.37}
|I^+_2(\lambda)|:=|\int_{-\infty}^\infty e^{i\lambda
t}Tr\{f(u^1(t)-u^2(t))\}\rho_2^+(t)dt|\leq
C_N(1+|\lambda|)^{-N}e^{a^+(-\Im \lambda)};
\end{equation}
if supported on $\mathbb{R}^-$, we choose $\Im\lambda<0$, the
lower half complex plane, for
\begin{equation}\label{2.38}
|I^-_2(\lambda)|:=|\int_{-\infty}^{\infty}e^{i\lambda
t}Tr\{f(u^1(t)-u^2(t))\}\rho_2^-(t)dt|\leq
C_N(1+|\lambda|)^{-N}e^{a^{-}(-\Im \lambda)},
\end{equation}
where $a^{\pm}$ is the supporting function of $\rho^{\pm}(t)$.
This form of Paley-Wiener's theorem appears in H\"{o}rmander's
book \cite{Hormander}. In this case,
\begin{equation}
|I_2(\lambda)|\leq C_N(1+|\lambda|)^{-N},\forall
N\in\mathbb{N},\mbox{ whenever }\lambda\in0i+\mathbb{R}.
\end{equation}
This is a rapidly decreasing term.

\par
For $I_3(\lambda)$, we see $\rho_3(t)$ is also an union of two
cutoff functions supported on $(-\infty,-\beta)$ and
$(\beta,\infty)$ respectively. By domain of dependence argument on
$u^k(t,x,x)$ along the geodesic toward $0\in\mathbb{R}^n$ hitting
the obstacle boundaries and back to $x$ such that
$\omega=-\theta$, we choose $\beta$ large such that
\begin{equation}
f(u^1(t,r\theta,r\theta)-u^2(t,r\theta,r\theta))\rho_3(t)\equiv0,\forall\theta\in\mathbb{S}^{n-1}.
\end{equation}
There are infinitely many geodesics starting at $x$ and back to
$x$. We consider only the one carries backscattering information.
Under starlike assumption, all such geodesics are transversal
reflections.

\par
Henceforth,
\begin{equation}
{\rm
Tr}f(u^1(t)-u^2(t))\rho_3(t)=\int_0^\infty\int_{\mathbb{S}^{n-1}}f(r\theta)(u^1(t,r\theta,r\theta)-u^2(t,r\theta,r\theta))\rho_3(t)r^{n-1}d\theta
dr\equiv0.
\end{equation}
Therefore,
\begin{equation}
I_3(\lambda)\equiv0.
\end{equation} $\Box$
\end{proof}
Accordingly,
\begin{corollary}
$I_1(\lambda)$ and $I_2(\lambda)$ are entire functions.
\end{corollary}
\begin{proof}
We see that $\rho_1(t)+\rho_2(t)
\in\mathbb{C}^\infty_0(\mathbb{R}).$ $\Box$
\end{proof}
\section{Proof of Theorem 1.1}
Let us define
\begin{equation}
\Phi_k(t):=\mathcal{F}_{\lambda\rightarrow
t}[|\lambda|^k],\label{phi}
\end{equation}
where one needs to replace $|\lambda|^k$ by its certain
regularization when $k\leq -1$. According to Ivrii
\cite{Ivrii,Ivrii2}, when $t\rightarrow 0^+$, we have
\begin{equation}
Tr\{f\cos
t\sqrt{\Delta_{\mathcal{O}}}\rho_1(t)\}=\sum_{j=0}^\infty
c_j\Phi_{n-j-1}(t),\label{3.16}
\end{equation}
where $c_j$'s are nonzero multiples of heat invariants
$a_{j/2}$'s. See Branson and Gilkey \cite{Branson}. In particular,
\begin{equation}\label{3.3}
c_0=\alpha_0\int_\mathcal{O} f(x)dx, \mbox{ where the constant
}\alpha_0\neq0.
\end{equation}

Let
\begin{equation}
D(\lambda):=\int_{-\infty}^\infty e^{i\lambda
t}Tr\{fu^1(t)\}\rho_1(t)dt-\int_{-\infty}^\infty e^{i\lambda
t}Tr\{fu^2(t)\}\rho_1(t)dt.
\end{equation}
By Lemma~\ref{11} and Proposition~\ref{2.3}, as a result of
analytic continuation,
\begin{eqnarray}\label{3.5}
D(\lambda)\mbox{ is rapidly decreasing on }0i+\mathbb{R}.
\end{eqnarray}
Using~(\ref{3.16}), on the other hand, as $\lambda\rightarrow
0i\hspace{-2pt}\pm\infty$,
\begin{eqnarray}
D(\lambda)\rightarrow
\alpha_{0}(a_{0}(f,\mathcal{O}^1)-a_{0}(f,\mathcal{O}^2))|\lambda|^{n-1}
+\alpha_{0.5}(a_{0.5}(f,\mathcal{O}^1)-a_{0.5}(f,\mathcal{O}^2))|\lambda|^{n-2}+\cdots.\label{2.19}
\end{eqnarray}
Combing~(\ref{3.5}) and~(\ref{2.19}), we obtain
\begin{equation}
a_j(f,\mathcal{O}^1)=a_j(f,\mathcal{O}^2),\hspace{2pt}\forall
j\geq0,\mbox{ where
}f\in\mathcal{C}^\infty_0(\Omega^1\cap\Omega^2).\label{43}
\end{equation}
In particular, we have identical localized relative volume
\begin{equation}
a_0(f,\mathcal{O}^1)=a_0(f,\mathcal{O}^2),\label{47}
\end{equation}
where we choose that $f=f(x)=f(r\omega)=f(\omega)$, where
$\omega\in\mathbb{S}^{n-1}$. It suffices to show the obstacle can
be shaped by angular argument. By our starlike assumption, we have
\begin{equation}
r^k(\omega):=\sup\{v|v\omega\in\mathcal{O}^k\}\end{equation} as
the radial function of $\mathcal{O}^k$ in the direction of
$\omega\in\mathbb{S}^{n-1}$. According to \cite{Groemer} and
\cite[equation(2.4)]{Koldobsky}, we have for starlike sets
\begin{equation}
{\rm
Volume}(\mathcal{O}^k)=\int_{\mathbb{S}^{n-1}}(r^k)^n(\omega)d\omega.
\end{equation}
 Hence,~(\ref{3.3}) and~(\ref{47}) becomes
\begin{equation}
\int_{\mathbb{S}^{n-1}}(r^1)^n(\omega)f(\omega)d\omega=\int_{\mathbb{S}^{n-1}}(r^2)^n(\omega)f(\omega)d\omega,
\forall f\in\mathcal{C}^\infty_0(\mathbb{S}^{n-1}).
\end{equation}Therefore, we have $r_1(\omega)=r_2(\omega)$. Theorem is proved.


\begin{thebibliography}{widest-label}
\bibitem{Bers}L. Bers, F. John and M. Schechter, Partial
differential equation, Lectures in applied mathematics, V.3A,
1964, Johan Wiley and sons.
\bibitem{Branson}T.P. Branson and P.B. Gilkey,
The asymptotics of the Laplacian on a manifold with boundary,
 Comm. Partial Differential Equations, 15(1990), no. 2, 245-272.


\bibitem{Balslev}E. Balslev, Absence of positive eigenvalues of Schr\"{o}dinger
operators, Archive for rational mechanics and analysis, V.59,
Number 4, 343-357(1975).


\bibitem{Colton and Sleeman}D. Colton and B.D. Sleeman, Uniqueness
theorems for the inverse problem of acoustic scattering, IMA
Jounal of Applied Mathematics, 31(1983), 253-259.
\bibitem{Groemer}H. Groemer, Geometric applications of Fourier series and spherical
harmonics, Encyclopedia of mathematics and its applications, v.
61,  Cambridge University Press, New York, 1996.





\bibitem{Hormander}L. H\"{o}rmander, The analysis of linear
partial differential operators I, Springer-Verlag,
Berlin-Heidelberg, 1990.
\bibitem{Isakov}V. Isakov, New stability results for soft obstacles
in inverse scattering, Inverse Problem, 9(1993), 535-543.
\bibitem{Isakov2}V. Isakov, Uniqueness and stability in
multi-dimenstional inverse problem, Inverse Problems, 9(1993),
579-621.
\bibitem{Ivrii}V. Ivrii, Precise spectral asymptotics for elliptic
operators acting in fiberings over manifolds with boundary,
Lecture notes in mathematics, V.1100, Springer-Verlag, Berlin
Heidelberg New York Tokyo, 1984.
\bibitem{Ivrii2}V. Ivrii, Second term of the spectral asymptotic
expansion of the Laplace-Beltrami operator on manifolds with
boundary, Funktsioal'nyi Analiz i Ego Prilozheniya, 14(1980),
No.2, 25-34.
\bibitem{Ivrii3}V. Ivrii, Exact spectral asymptotics for the
Laplace-Beltrami operator in the case of general elleptic boundary
conditions, Funksional'nyi Analiz i Ego Prilozheniya, V 15, No.1,
74-75(1981).


\bibitem{Koldobsky}A. Koldobsky, Fourier analysis in convex geometry, Mathematical Surveys
and Monographs, 116, American Mathematical Society, Providence,
RI, 2005.

\bibitem{Lax and Phillips}P.D. Lax and R.S. Phillips, Scattering theory, New York, Acdemic press, 1989.
\bibitem{Lax and Phillips2}P.D. Lax and R.S. Phillips, A logarithmic bound on the location of
the poles of the scattering matrix, Arch. Rational Mech. Anal, 40,
268-280(1971).

\bibitem{Melrose}R.B Melrose, Geometric scattering theory, Cambridge university
press, 1995.
\bibitem{Petkov}V. Petkov and G. Popov, Asymptotic
behavior of the scattering phase for non-trapping obstacles, Ann.
Inst. Fourier, Grenoble, 32(3), 111-149(1982).

\bibitem{Reed and Simon}Reed and Simon, Methods of Modern
Mathematical Physics, v.1 and v.2 , Academic press, new York,
1975.

\bibitem{Shenk}N. Shenk and D. Thoe, Resonant states and poles of
the scattering matrix for perturbations of $-\Delta$, Journal of
mathematical analysis and applications, 37, 467-491(1972).
\bibitem{Shenk2}N. Shenk and D. Thoe, Outgoing solutions of
$(-\Delta+q-k^2)u=f$ in an exterior domain, Journal of
mathematical analysis and applicaitons, 31, 81-116(1970).
\bibitem{Sjostand and Zworski}J. Sj\"{o}strand and M. Zworski,
Complex scaling method and the distribution of scattering poles,
J. Amer. Math. Soc. 4(1991), 729-769.
 \bibitem{Titch}E.C. Titchmarsh, "The Theory of Functions", Oxford
University Press, 2nd ed.


\bibitem{Zworski3}M. Zworski, Poisson formula for resonances in
even dimension, Asian J. Math, 2(3)(1998), 615-624.
\end{thebibliography}
\end{document}